\documentclass[12pt]{article}

\usepackage{amssymb}
\usepackage{amsmath}

 \newtheorem{thm}{Theorem}[section]
 \newtheorem{prop}{Proposition}[section]

 \newtheorem{rem}{Remark}[section]

\renewcommand{\leq}{\leqslant}
\renewcommand{\geq}{\geqslant}
\renewcommand{\setminus}{\smallsetminus}

\newcommand{\rbx}{\hfill{\rule{1ex}{1ex}}}

\newcommand{\ind}{\mathrm{ind}\,}

 \newcommand{\vp}{\varphi}

\newcommand{\coker}{\mbox{\rm coker}\,}

\newcommand{\im}{\mathrm{im}\,}

\newcommand{\nn}{\nonumber}

\newcommand{\cF}{\mathcal{F}}

\newcommand{\cL}{\mathcal{L}}

\newcommand{\cR}{\mathcal{R}}

\newcommand{\fB}{{\mathfrak B}}

\newcommand{\sC}{{\mathbb C}}

\newcommand{\sN}{{\mathbb N}}

\newcommand{\sR}{{\mathbb R}}

\newcommand{\sZ}{{\mathbb Z}}

\begin{document}

\vspace*{10mm}

\begin{center}
{\Large\textbf{Kernels of Wiener-Hopf plus Hankel operators}}

\vspace{3mm}

{\Large\textbf{with matching
generating functions}}
\end{center}

\vspace{5mm}

\begin{center}

\textbf{Victor D. Didenko and Bernd Silbermann}

\vspace{2mm}


Universiti Brunei Darussalam, Bandar Seri Begawan, BE1410  Brunei;
diviol@gmail.com

Technische Universit{\"a}t Chemnitz, Fakult{\"a}t f\"ur Mathematik,
09107 Chemnitz, Germany; silbermn@mathematik.tu-chemnitz.de

 \end{center}

\vspace{5mm}

\textit{Dedicated to Roland Duduchava on the occasion of his seventieth
birthday}

\begin{abstract}
Considered are Wiener--Hopf plus Hankel operators
$W(a)+H(b):L^p(\sR^+)\to L^p(\sR^+)$ with generating functions  $a$
and $b$ from a subalgebra of $L^\infty(\sR)$ containing almost
periodic functions and Fourier images of $L^1(\sR)$-functions. If
the generating functions  $a$ and $b$ satisfy the matching condition
 \begin{equation*}
a(t) a(-t)=b(t) b(-t),\quad t\in\sR,
 \end{equation*}
an explicit description for the kernels and cokernels of the
operators mentioned is given.
 \end{abstract}

  \vspace{6mm}


\vspace{2mm}


\section{Introduction\label{s1}}

Let $\sR^-$ and $\sR^+$ be, respectively, the subsets of all
negative and all positive real numbers, and let $\chi_E$ refer to
the characteristic function of the subset $E$ of the set of real
numbers $\sR$, i.e.
 $$
\chi_E(t):=\left \{
\begin{array}{ll}
  1 & \text{ if } t\in E,  \\
  0 & \text{ if } t\in \sR\setminus E. \\
\end{array}
 \right.
 $$
By $L^p(\sR^+):=\chi_{\sR^+} \, L^p(\sR)$ and
$L^p(\sR^-):=\chi_{\sR^-}\, L^p(\sR)$ we denote the subspaces of
$L^p(\sR)$, $1\leq p \leq\infty$ which contain all functions
vanishing on $\sR^-$ and $\sR^+$, cor\-respondingly.

Consider the set $G$ of functions defined on the real line $\sR$ and
having the form
  \begin{equation}\label{cst1}
  a(t)=\sum_{j=-\infty}^\infty a_j e^{i\delta_j t}
 +\int_{-\infty}^\infty k(s) e^{its} \, ds, \quad -\infty<t<\infty,
 \end{equation}
where $\delta_j$ are pairwise distinct real numbers and
  \begin{equation*}
  \sum_{j=-\infty}^\infty |a_j|<\infty, \quad
 \int_{-\infty}^\infty |k(s)| \, ds <\infty.
 \end{equation*}
Each element $a$ of $G$ generates three operators
$W^0(a):L^p(\sR)\to L^p(\sR)$ and $W(a), H(a):L^p(\sR^+)\to
L^p(\sR^+)$,
 \begin{align*}
 (W^0(a)f)(t)&:=\sum_{j=-\infty}^\infty a_j f(t-\delta_j)
 +\int_{-\infty}^\infty k(t-s) f(s)\,ds, \\
  W(a)&:=PW^0(a) , \nn \\
  H(a)&:= PW^0(a)QJ, \nn
\end{align*}
where $P: f\to \chi_{\sR^+} f$ and $Q:=I-P$ are the projections on
the subspaces $L^p(\sR^+)$ and $L^p(\sR^-)$, correspondingly, and
the operator $J:L^p(\sR)\to L^p(\sR)$ is defined by  $J\vp :=
\widetilde{\vp}$ with $\widetilde{\vp}(t):=\vp(-t)$. Note that
$W^0(a), W(a)$ and $H(a)$ are bounded linear operators on the
corresponding spaces. The function $a$ is called the generating
function or the symbol for each of the operators $W^0(a)$, $W(a)$
and $H(a)$. Wiener-Hopf and Hankel operators are closely connected.
Thus for any $a,b\in G$, one has
\begin{equation}\label{cst4}
 \begin{aligned}
    W(ab)&=W(a)W(b)+H(a) H(\widetilde{b}),\\
    H(ab)&= W(a)H(b)+H(a)W(\widetilde{b}).
\end{aligned}
 \end{equation}
The Fredholm theory for the operators $W^0(a)$, $a\in G$ is
relatively simple. An operator $W^0(a)$ is semi-Fredholm if and only
if $a$ is invertible in $G$. The study of the operators $W(a)$ is
more involved. Nevertheless, for various classes of generating
functions $a$, Wiener-Hopf operators $W(a)$ are well studied  (see,
for example, \cite{BS:2006,BKS:2002,CD:1969,Du:1973,
Du:1977,Du:1979,GF1974}). In particular, Fredholm properties of such
operators are known and a description of the kernel is available. On
the other hand, Wiener-Hopf plus Hankel operators, i.e. the
operators of the form
 $$
B=B(a,b)= W(a)+H(b)
 $$
remains less studied. Fredholm properties of such operators can be
derived by reducing the initial operator to a Wiener-Hopf operator
with a matrix symbol, and there is a number of works where this idea
is successfully implemented
\cite{MR2400132,CN:2009,CS2010d,CS:2011}. However, these works
mainly deal with generating functions $a$ and $b$ satisfying the
condition $a=b$ and consider the operators acting in an $L^2$-space.
If $a\neq b$, then a rarely verifiable assumption about special
matrix factorization is used. A different approach to the study of
the operators of the form $I+H(b)$ has been employed in
\cite{KS2000, Karapetiants2001}, where the essential spectrum and
the index of such operators have been found. On the other hand, no
information is available about the kernel elements of the operators
$W(a)+H(b)$ even in the above mentioned cases $a=1$ and $a=b$. The
goal of this work is to present an efficient description of the
space $\ker B(a,b)$ when the generating functions $a$ and $b$ belong
to the Banach algebra $G$ and satisfy a specific algebraic relation.
Point out that our approach does not involve factorization of any
matrix function but only the one of scalar functions.

Let $a,b\in L^\infty(\sR)$. We say that the duo $(a,b)$ is a
matching pair if
\begin{equation}\label{cst9}
   a \widetilde{a}= b \widetilde{b},
\end{equation}
where $\widetilde{a}:=a(-t)$. The relation \eqref{cst9} is called
matching condition. In the following we always assume that $a$, and
therefore $b$, is invertible in $G$. For each matching pair $(a,b)$,
consider the pair $(c,d)$ with
  $$
 c:=\widetilde{b}\widetilde{a}^{-1}, \quad d:= b \widetilde{a}^{-1}.
  $$
It is easily seen that $(c,d)$  is also a matching pair. This pair
is called the subordinated pair for $(a,b)$ or just the subordinated
pair. The elements $c$ and $d$ of the subordinated pair possess a
specific property, namely
  $$
c\tilde{c}=1, \quad d\widetilde{d}=1.
 $$
Throughout this paper any function $g\in G$ satisfying the condition
 \begin{equation*}
g \widetilde{g}=1,
 \end{equation*}
is called  matching function. Note that the matching functions $c$
and $d$ can be also expressed as
\begin{equation*}
c=ab^{-1}, \quad d=  a \widetilde{b}^{-1}.
\end{equation*}
Further, if $(c,d)$ is the  subordinated pair for $(a,b)$, then
$(\overline{d},\overline{c})$ is the subordinated pair for the
matching pair $(\overline{a}, \overline{\widetilde{b}})$. Moreover,
if $p\in[1,\infty)$, then $\overline{a}$ and
$\overline{\widetilde{b}}$ are generating functions for the operator
adjoint to the Wiener-Hopf plus Hankel operator
$W(a)+H(b):L^p(\sR^+)\to L^p(\sR^+)$, i.e.
  \begin{equation}\label{cst10}
(W(a)+H(b))^*=W(\overline{a})+H(\overline{\widetilde{b}}).
 \end{equation}
The Wiener-Hopf operators with matching generating symbols possess a
number of remarkable properties. In particular, the kernels of such
operators can be structured in a special way and this
structurization can be used in the description of the kernels of
Wiener-Hopf plus Hankel operators. More precisely, let $g$ be a
matching function and let $\mathbf{P}(g)$ be the operator defined on
the kernel $\ker W(g)$ by
 \begin{equation}\label{Pro1}
\mathbf{P}(g):=J QW^0(g)P\left |_{\ker W(g)}\right ..
 \end{equation}
One can easily check that $\mathbf{P}(g)$ maps $\ker W(g)$ into
$\ker W(g)$ and $\mathbf{P}^2(g)=I$ (see \cite{DS:2014b} for more
details). Therefore, the operators
 \begin{equation}\label{Pro2}
\mathbf{P}^{-}(g):=(1/2)(I- \mathbf{P}(g)), \quad
\mathbf{P}^{+}(g):=(1/2)(I+ \mathbf{P}(g)),
 \end{equation}
considered on the space $\ker W(g)$, are complementary projections
generating a decomposition of $\ker W(g)$, i.e.
 \begin{equation*}
   \ker W(g)=\im \mathbf{P}^{-}(g)\dotplus\im\mathbf{P}^{+}(g).\\
   \end{equation*}
Consider now the Wiener-Hopf plus Hankel operators $W(a)+H(b)$,
generating functions of which constitute a matching pair. In this
case the elements of the subordinated pair $(c,d)$ are matching
functions. Assume that the operator $W(c)$ is right-invertible and
let $W_r^{-1}$ be a right inverse for $W(c)$. By $\vp_{\pm}$ we
denote the operators defined on the kernel of the operator $W(d)$ by
  \begin{equation}\label{eqphi}
   2\vp_{\pm}(s):=W_r^{-1}(c)W(\widetilde{a}^{-1}) s \mp J Q W^0(c)P
   W_r^{-1}(c)W(\widetilde{a}^{-1})s \pm J Q W^0(\widetilde{a}^{-1})
   s,
\end{equation}
where $\widetilde{a}^{-1}=a^{-1}(-t)$. It was shown in
\cite{DS:2014b} that for any $s\in \ker W(d)$ one has
$\vp_{\pm}(s)\in \ker (W(a)\pm H(b))$, and the operators $\vp_+$ and
$\vp_-$ are injections on the spaces $\im \mathbf{P}^+(d)$ and $\im
\mathbf{P}^-(d)$, respectively. Moreover, the following result is
true.
   \begin{prop}[{see \cite[Proposition 2.3]{DS:2014b}}]\label{p2}
 Let $(c,d)$ be the subordinated pair for a matching pair $(a,b)\in G \times G$.
 If the operator $W(c)$ is right-invertible, then
 \begin{equation}\label{cst14}
\begin{aligned}
   \ker(W(a)+H(b))& =\vp_{+}(\im \mathbf{P}^+(d)) \dotplus\im \mathbf{P}^-(c), \\
\ker(W(a)-H(b))& =\vp_{-}(\im \mathbf{P}^-(d)) \dotplus\im
\mathbf{P}^+(c).
\end{aligned}
\end{equation}
 \end{prop}

Thus to  describe the kernels of the Wiener-Hopf plus/minus Hankel
operators, one needs to find an efficient description of the images
of the projections $\mathbf{P}^{\pm}(c)$ and $\mathbf{P}^{\pm}(d)$.
Notice that the above statements do not depend on $p$.

This paper is organized as follows. In Section \ref{s2} we present a
decomposition of the kernel of $W(g)$ with a generating matching
function $g$. These results are used in Section \ref{s3} in order to
derive an efficient description of the kernels $\ker(W(a)\pm H(b))$,
$p\in [1,\infty]$ and the cokernels $\coker (W(a)\pm H(b))$, $p\in
[1,\infty)$. Similar results for Toeplitz plus Hankel operators have
been obtained in \cite{DS:2014a, DS:2014}, and generalized Toeplitz
plus Hankel operators are considered in \cite{DS:2013i}. However,
all the relevant operators in \cite{DS:2014a, DS:2014, DS:2013i} are
Fredholm. On the other hand, the really new feature of the present
study is the consideration of situations where the operators $W(c)$
and $W(d)$ can have infinite-dimensional kernels or co-kernels.

\section{Kernels of Wiener-Hopf operators with a matching
generating function.\label{s2}}

Our aim now is to describe the subspaces $\im
\mathbf{P}^{\pm}(g)\subset \ker W(g)$. For, let us recall certain
results of Fredholm theory for Wiener-Hopf operators with generating
functions from the Banach algebra $G$. As we know, any element $a\in
G$ can be represented in the form $a=b+k$, where $b$ belongs to the
algebra $AP_w$ of all almost periodic functions with absolutely
convergent Fourier series and $k$ is in the algebra $\cL_0$ of all
Fourier transforms of functions from $L^1(\sR)$. If $a=b+k$, $b\in
AP_w, k\in \cL_0$ is an invertible element of $G$, then $b$ is
invertible in $AP_w$ and one can define the numbers $\nu=\nu(a)$ and
$n=n(a)$ by
  \begin{equation*}
   \nu(a):=\lim_{l\to\infty}\frac{1}{2l} [\arg b(t)]_{-l}^l,\quad
     n(a):=\frac{1}{2\pi} [\arg
     (1+b^{-1}(t)k(t)]_{t=-\infty}^\infty.
 \end{equation*}
Recall that $a\in G$ is invertible in $G$ if and only if
$\inf_{t\in\sR} |a(t)|>0$ and $\cL_0$ forms a closed two-sided ideal
in $G$.

  \begin{thm}[Gohberg/Feldman \cite{GF1974}]\label{t1}
Let $1\leq p\leq \infty$ and $g\in G$. The operator $W(g)$ is
one-sided invertible in the space $L^p(\sR^+)$ if and only if $g$ is
invertible in $G$. Further, if $g\in G$ is invertible in $G$, then
the following assertions are true:
 \sloppy
\begin{enumerate}

        \item If $\nu(g)<0$, then the operator $W(g)$ is invertible
    from the right and $\dim \ker W(g)=\infty$.
        \item If $\nu(g)=0$ and $n(g)\geq 0$ ($\nu(g)=0$ and $n(g)\leq 0$), then the operator $W(g)$ is invertible
    from the left (from the right) and
      \begin{equation*}
      \dim \coker W(g)=n(g) \quad (\dim \ker W(g)=-n(g)).
      \end{equation*}
      \item  If $\nu(g)>0$, then the operator $W(g)$ is invertible
    from the left and $\dim \coker W(g)=\infty$.
      \item If $g\in G$ is not invertible in $G$, then $W(g)$ is not
      a semi-Fredholm operator.
\end{enumerate}
  \end{thm}
The proof of this theorem is based on the fact that every invertible
function $a\in G$ admits a factorization of the form
  \begin{equation}\label{EqFactor}
  g(t)=g_-(t) e^{i\nu t} \left (\frac{t-i}{t+i} \right )^n g_+(t), \quad
  -\infty<t<\infty,
 \end{equation}
where $g_{+}^{\pm1}\in G^+$, $g_{-}^{\pm1}\in G^-$,  $\nu=\nu(g)$
and $n=n(g)$. Recall that $G^+(G^-)$ is defined as follows:
$G^+(G^-)$ consists of all functions \eqref{cst1} such that all
indices $\delta_j$ are non-negative (non-positive) and the function
$k$ vanishes on the negative (positive) semi-axis. It is clear that
functions from $G^+$ and $G^-$ admit holomorphic extensions to the
upper and to the lower half-plane, correspondingly, and the
intersection of the algebras $G^+$ and $G^-$ consists of constant
functions only. Note that under the condition $g_-(0)=1$, the
factorization \eqref{EqFactor} is unique. Moreover, for $a\in G^-, b
\in G$ and $c\in G^+$, the first identity from \eqref{cst4} leads to
the relation
 $$
W(abc)=W(a) W(b) W(c).
 $$
Combined with the factorization \eqref{EqFactor}, this relation
leads to the following representation of the operator $W(g)$,
 $$
W(g)= W(g_-) W\left ( e^{i\nu t} \left ( \frac{t-i}{t+i}\right )^n
\right ) W(g_+).
 $$
Therefore, theory of the Wiener-Hopf operators $W(g)$ with
invertible symbol $g$ is based on the study of the middle factor of
this factorization (see \cite[Chapter VII]{GF1974}). Thus the
operator $W(a)$ has a kernel containing non-zero elements in the two
cases--viz. if $\nu<0$, then $\dim\ker W(g)=\infty$, or if $\nu=0$
and $n<0$, then $\dim\ker W(g)=|n|$. In what follows we consider all
possible situations separately. Let us note that $\ker W(a)$ do not
depend on $p$.

Assume that $g$ is a matching function.  Then, as was pointed out in
\cite{DS:2014b}, the factorization \eqref{EqFactor} comes down to
the following one
  \begin{equation}\label{Eq1}
  g(t)= \boldsymbol\sigma(g)\, \widetilde{g}_+^{-1}(t) e^{i\nu t} \left (\frac{t-i}{t+i} \right )^n
  g_+(t)
 \end{equation}
where $\boldsymbol\sigma(g)=(-1)^n g(0)$,
$\widetilde{g}_+^{\pm1}(t)\in G^-$ and
$g_-(t)=\boldsymbol\sigma(g)\,\widetilde{g}_+^{-1}(t)$. In passing
note that $\boldsymbol\sigma(g)=\pm1$.


Our goal now is to describe the projections $\mathbf{P}^{\pm}(g)$
from \eqref{Pro2}. Let us start with the case where the parameters
$\nu$ and $n$ in the factorization \eqref{Eq1} satisfy the relations
$\nu=0$, $n<0$. It is known \cite{GF1974} that in this case
  \begin{equation*}
\ker W(g)=\left \{ W(g_+^{-1}) \left ( \sum_{j=0}^{|n|-1}c_j t^j
 e^{-t} \right ): c_j\in \sC \right
  \}.
 \end{equation*}
Thus the functions $W(g_+^{-1})t^j e^{-t}$, $j=0,1,\cdots,|n|-1$
form a basis in $\ker W(g)$. However, the space $\ker W(g)$ has
another basis, namely,
 \begin{equation}\label{Eq1.5}
\{ W(g_+^{-1}) \psi_j(t): \quad j=0,1,
 \cdots, |n|-1\},
 \end{equation}
 where
    \begin{align*}
    \psi_j(t)&:= \left \{
\begin{array}{ll}
\sqrt{2} e^{-t} \Lambda_j(2t),& \text{ if } t>0,\\
 0, & \text{ if } t<0,\\
 \end{array}
     \right .,& \quad j=0,1,\cdots\,,\phantom{--}
     \end{align*}
and $\Lambda_j$ are the normalized Laguerre polynomials. Moreover,
for $j=-1,-2, \ldots,$ one can define the functions $\psi_j$ by
 \begin{align}
    \psi_j(t) &:=-\psi_{-j-1}(-t), &\quad j=-1,-2, \cdots \,, \label{Eq3}
\end{align}
The functions $\psi_j$, $j\in \sZ$ can be also expressed in the form
  \begin{equation}\label{Eq4}\begin{aligned}
    \psi_j(t)&=(U^j \psi_0)(t), \quad j=\pm1,\pm2, \cdots \,,\\
    \psi_0(t)&= \left \{
    \begin{array}{ll}
\sqrt{2} e^{-t}, & \text{ if } t>0\\
 0, & \text{ if } t<0,
    \end{array}
     \right .\, .
    \end{aligned}
 \end{equation}
where $U:=W^0((\lambda-i)/(\lambda+i))$. Note that the operators
$U^j,j\in \sZ$ are unitary operators on $L^2(\sR)$. Thus, the
functions $\psi_j, j\in\sZ$  form an orthonormal basis on this
space. Indeed, it is shown in \cite[Chapter 3, \S 3.2]{GF1974} that
for $j>0$, one has
 $$
(U^j \psi_0)(t)=\psi_j(t),
 $$
and applying \eqref{Eq3} one gets the result. Note that the relation
\eqref{Eq3} can be obtained by using the Fourier transform. Indeed,
let us recall the formula
 $$
(\cF \psi_n)(\lambda)=\int_0^\infty \psi_n (t)\,e^{i\lambda
t}\,dt=\int_0^\infty U^n \psi_0 (t)\,e^{i\lambda t}\,dt= \left
(\frac{\lambda-i}{\lambda+i}\right )^n \frac{i\sqrt{2}}{\lambda+i},
\quad n\in\sZ_+,
 $$
where $\cF$ is the Fourier transform \cite{GF1974} and $\sZ_+$
refers to the set of all non-negative integers. Consider the
operator $J:L^p(\sR)\to L^p(\sR)$ defined by $(Jf)(t)=f(-t)$. If
$n\in \sN$, then one has
\begin{align}\label{psi1}
 \cF(-J\psi_{n-1})(\lambda)=-J\cF(\psi_{n-1})(\lambda)
\nn \\
=-J\left ( \left (\frac{\lambda-i}{\lambda+i} \right
)^{n-1}\frac{i\sqrt{2}}{\lambda+i} \right ) &=-\left ( \left
(\frac{-\lambda-i}{-\lambda+i} \right
)^{n-1}\frac{i\sqrt{2}}{-\lambda+i} \right )\nn \\
&= \left (\frac{\lambda+i}{\lambda-i} \right
)^{n-1}\frac{i\sqrt{2}}{\lambda-i}.
\end{align}
On the other hand, if $n\leq-1$, then
 \begin{align}\label{psi2}
(\cF\psi_n)(\lambda)& =\int_0^\infty \psi_n(t)e^{i\lambda t}\, dt
 =\int_0^\infty U^n \psi_0(t)e^{i\lambda t}\, dt \nn\\
  &=
\left (\frac{\lambda+i}{\lambda-i} \right
)^{|n|}\frac{i\sqrt{2}}{\lambda+i} =\left
(\frac{\lambda+i}{\lambda-i} \right
)^{|n|-1}\frac{i\sqrt{2}}{\lambda-i}.
 \end{align}
Comparing \eqref{psi1} and \eqref{psi2}, one obtains that
 $$
\cF(\psi_n (t))=\cF(-\psi_{|n|-1}(-t))
 $$
and one has to use the injectivity of the Fourier transform to
complete the proof.

Let $g$ be a matching function. In order to describe the
corresponding projections $\mathbf{P}^{\pm}(g)$ of
\eqref{Pro1}-\eqref{Pro2}, we will study how the operator
$\mathbf{P}(g)$ interacts with the basis elements \eqref{Eq1.5}.
Thus
\begin{align*}
   \mathbf{P}(g) W(g_+^{-1}) \psi_j (t)& =JQW^0(g)P
   W(g_+^{-1})\psi_j(t)\\
    &=JQW^0 \left(\boldsymbol\sigma(g)\, \widetilde{g}_+^{-1}  \left (\frac{t-i}{t+i} \right )^{-|n|}
  g_+  \right ) W(g_+^{-1}) \psi_j  \\
&=\boldsymbol\sigma(g)\,JQW^0(\widetilde{g}_+^{-1} ) W^0\left (\left
(\frac{t-i}{t+i} \right )^{-|n|}\right ) \psi_j  .
\end{align*}
Considering the elements $W^0\left (\left ((t-i)/(t+i) \right
)^{-|n|}\right ) \psi_j$, $j=0,1,\cdots, |n|-1$ and using relations
\eqref{Eq3} and \eqref{Eq4}, we get
 \begin{align*}
W^0\left (\left (\frac{t-i}{t+i} \right )^{-|n|}\right ) \psi_j
 &=
  W^0\left (\left (\frac{t-i}{t+i} \right )^{-|n|}\right )W^0\left (\left (\frac{t-i}{t+i} \right
  )^j \right ) \psi_0\\
   &=W^0\left (\left (\frac{t-i}{t+i} \right
  )^{-|n|+j} \right ) \psi_0 =\psi_{-|n|+j}=-J\psi_{|n|-j-1}.
  \end{align*}
Hence,
 $$
\boldsymbol\sigma(g)\,JQW^0(\widetilde{g}_+^{-1} ) W^0\left (\left
(\frac{t-i}{t+i} \right )^{-|n|}\right ) \psi_j
=-\boldsymbol\sigma(g)\,PW^0(g_+^{-1} )\psi_{|n|-j-1}.
 $$
Now one can proceed similarly to \cite[Section 5]{DS:2014} and
obtain the following result.
 \begin{thm}\label{thm1}
Let $g\in G$ be a matching function such that the operator
$W(g):L_p(\sR_+)\to L_p(\sR_+)$ is Fredholm and $n:=\ind W(g)>0$. If
 $$
g(t)=  g_-(t) \left ( \frac{t-i}{t+i} \right)^{-n}g_+(t)=
\boldsymbol\sigma(g)\, \widetilde{g}_+^{-1}(t) \left
(\frac{t-i}{t+i} \right )^{-n} g_+ (t)\,  , \quad g_-(0)=1,
 $$
is the related Wiener-Hopf factorization of the function $g$, then
the following systems $\fB_{\pm}(g)$ of functions $W(g_+^{-1})
\psi_{j}$ form bases in the spaces $\im \mathbf{P}^{\pm}(g)$:
 \begin{enumerate}
\item If $n=2m$, $m\in\sN$, then
 $$
\fB_{\pm}(g)=\{W(g_+^{-1}) \left ( \psi_{m-k-1}\mp
\boldsymbol\sigma(g)\psi_{m+k}\right ): k=0,1,\cdots, m-1\},
 $$
 and
  $$
\dim\im \mathbf{P}^{\pm}(g)=m.
  $$
\item If $n=2m+1$, $m\in\sZ_+$, then
 $$
\fB_{\pm}(g)=\{W(g_+^{-1})\left ( \psi_{m+k}\mp
\boldsymbol\sigma(g)\psi_{m-k}\right ): k=0,1,\cdots, m\},
 $$
 and
  $$
\dim\im \mathbf{P}^{\pm}(g)=m+ \frac{1\mp\boldsymbol\sigma(g)}{2}.
  $$
 \end{enumerate}
 \end{thm}

 \begin{rem}\label{rem1}
It is worth mentioning that the zero element belongs to the one of
the sets $\{W(g_+^{-1})(\psi_{m+k}- \boldsymbol\sigma
(g)\psi_{m-k}):k=0,1,\cdots, m
 \}$ or
$\{W(g_+^{-1})(\psi_{m+k}+ \boldsymbol\sigma
(g)\psi_{m-k}):k=0,1,\cdots, m
 \}$ only. Namely, for $k=0$ one of the terms
$\psi_m(1\pm\boldsymbol\sigma(g))$ is equal to zero.
 \end{rem}

Consider now the case $\nu<0$ and $n=0$.  Then
  \begin{equation*}
\ker W(g)= \left \{ W(g_+^{-1})f: f\in L^p(\sR^+) \text{ and }
f(t)=0 \text{ for } t>|\nu|\right \},
 \end{equation*}
(see \cite[Chapter VII, \S 2.4]{GF1974}).

 \begin{thm}\label{thm2}
Let $g\in G$ be a matching function such that the function $g$
possesses the Wiener-Hopf factorization
 $$
g(t)= g_-(t) e^{i\nu t} g_+(t)=
\boldsymbol\sigma(g)\,\widetilde{g}_+^{-1}(t) e^{i\nu t}\,g_+ (t) ,
\quad \nu<0 \text{ and } g_-(0)=1,
 $$
and let $h\in \ker W(g)$, that is $h=W(g_+^{-1})f$ with an $f\in
L^p(\sR^+)$ such that $f(t)=0$ for $t>|\nu|$. Then
 $$
JQW^0(g)Ph=\boldsymbol\sigma(g)\,W(g_+^{-1}) \cR_{|\nu|} f,
 $$
 where
 \begin{equation}\label{Eq4.5}
 (\cR_{|\nu|})(t)=\left \{
 \begin{array}{cc}
   f(|\nu|-t), & \text{ if } 0<t< |\nu|\\
   0 & \text{ if } t> |\nu| \\
 \end{array}
 \right .,
 \end{equation}
and
 $$
\mathbf{P}^{\pm}(g)h= \frac{h \pm \boldsymbol\sigma(g)\,W(g_+^{-1})
\cR_{|\nu|}\, f}{2}.
 $$
 \end{thm}
The prof of this result runs similarly to the proof of Theorem
\ref{thm4} below where a more general factorization of the
corresponding matching function $g$ has to be used.

Next we consider the situation $\nu<0$ and $n<0$. In this case the
function $g$ admits the Wiener-Hopf factorization of the form
  \begin{equation}\label{Eq5}
g=\boldsymbol\sigma(g)\,\widetilde{g}_+^{-1} e^{i\nu t}\left
(\frac{t-i}{t+i} \right )^{n} g_+ \,. 
 \end{equation}
As is shown in \cite[Chapter VII]{GF1974}, the kernel of the
operator $W(g)$ is the direct sum of the kernels of the operators
$W(g_+((t-i)/(t+i))^n)$ and $W(g_+ e^{i\nu t})$. Thus
 $$
\ker W(g)= \ker W\left (g_+\left (\frac{t-i}{t+i}\right )^n \right
)\dotplus \ker W\left (g_+ e^{i\nu t}\right ).
 $$
Therefore, in order to characterize the projections
$\mathbf{P}^{\pm}(g):\ker W(g)\to \ker W(g)$, one can describe their
action on the subspaces $\ker W(g_+((t-i)/(t+i))^n)$ and $ \ker
W\left (g_+ e^{i\nu t}\right )$ separately. To this aim, let us use
the following representations of the function $g$:
 \begin{align*}
g=e^{i\nu t} g_1, \quad g_1:=\boldsymbol\sigma
(g)\widetilde{g}_+^{-1} \left ( \frac{t-i}{t+i} \right )^n g_+,\\
  g=\left ( \frac{t-i}{t+i} \right )^n g_2, \quad g_2:=\boldsymbol\sigma
(g)\widetilde{g}_+^{-1} e^{i\nu t} g_+. \label{Eqg2}
\end{align*}
Moreover, observe that $JQW^0(g)P=H(\widetilde{g})$.

 \begin{thm}\label{thm3}
Assume that $g$ is a matching function of the form \eqref{Eq5}.
 \begin{enumerate}
    \item  If $h\in \ker W(g_+((t-i)/(t+i))^n)$, then
     $$
 \mathbf{P}^{\pm}(g)h= \frac12 \left  [I\pm \left (  W(e^{i|\nu|
t})(\mathbf{P}^+(g_1)-\mathbf{P}^-(g_1))\right ) \right ] h.
     $$

    \item  If $h\in \ker W(g_+e^{i\nu t})$, then
     $$
\mathbf{P}^{\pm}(g)h= \frac12\left  [I\pm \left ( W\left (\left (
\frac{t-i}{t+i} \right
)^{|n|}\right)(\mathbf{P}^+(g_2)-\mathbf{P}^-(g_2))\right ) \right ]
h.
     $$
 \end{enumerate}
 \end{thm}

\textbf{Proof.}
Let us start with assertion (i). Using  \eqref{cst4} we obtain
 $$
JQW^0(g)P=PW^0(\widetilde{g})QJ=H(\widetilde{g})= W \left
(\widetilde{e^{i\nu t}}\right )H(\widetilde{g}_1)+H\left
(\widetilde{e^{i\nu t}}\right ) W(g_1),
 $$
and the relation $W (g_1)h=0$ implies that
 $$
H(\widetilde{g})h=W \left (e^{i|\nu| t}\right )H(\widetilde{g}_1)h.
 $$
Therefore,
 \begin{align*}
\mathbf{P}^{\pm}(g)h&=\left [ \frac{I\pm H(\widetilde{g})}{2}\right
]h=\left [ \frac{I\pm W \left (e^{i|\nu| t}\right
)H(\widetilde{g}_1))}{2}\right ]h \\
&=\frac12 \left  [I\pm \left (  W(e^{i|\nu|
t})(\mathbf{P}^+(g_1)-\mathbf{P}^-(g_1))\right ) \right ] h,
\end{align*}
so the assertion (i) is proved.

The proof of assertion (ii) is similar to that of (i). It is based
on the formula
 $$
H(\widetilde{g})=W \left ( \widetilde{\left ( \frac{t-i}{t+i} \right
)^{|n|}}\right)  H(\widetilde{g}_2) + H\left ( \widetilde{\left (
\frac{t-i}{t+i} \right )^{|n|}}\right) W(g_2),
 $$
and is left to the reader.
\rbx

 \begin{rem}\label{rem2}
Recall that the projections $\mathbf{P}^{\pm}(g_1)$ and
$\mathbf{P}^{\pm}(g_2)$ acting, respectively, on the subspaces $\ker
W(g_1)$ and $\ker W(g_2)$ are described by Theorem \ref{thm1} and
Theorem \ref{thm2}. Besides,
 $$
\boldsymbol\sigma(g)=\boldsymbol\sigma(g_1)=\boldsymbol\sigma(g_2).
 $$
 \end{rem}

Finally, let us consider the case $\nu<0$ and $n>0$, i.e. now we
assume that the Wiener-Hopf factorization of the matching function
$g$ is
  \begin{equation}\label{Eq6}
g=\boldsymbol\sigma(g)\,\widetilde{g}_+^{-1} e^{i\nu t}\left
(\frac{t-i}{t+i} \right )^{n} g_+  \,. 
 \end{equation}
If this is the case, then according to \cite[Chapter VII]{GF1974}
the kernel of the operator $W(g)$ consists of the functions $h$
having the form
  \begin{equation}\label{Eq7}
  h=W(g_+^{-1})W\left ( \left ( \frac{t-i}{t+i}\right )^{-n} \right
  )\vp,
 \end{equation}
where $\vp\in L^p(\sR^+)$ is such that
 \begin{equation}\label{Eq8}
\vp(t)=0 \text{ for all } t>|\nu| \text{ and } \int_0^\infty
\vp(t)\, t^j e^{-t}\,dt=0, \quad j=0,1,\cdots, n-1.
 \end{equation}

 \begin{thm}\label{thm4}
Let $g\in G$ be a matching function such that the function $g$
possesses the Wiener-Hopf factorization \eqref{Eq6}. Assume that
$h\in \ker W(g)$. Then it can be represented in the form
\eqref{Eq7}--\eqref{Eq8} and
 $$
JQW^0(g)Ph=\boldsymbol\sigma(g)\,W(g_+^{-1}) \cR_{|\nu|} \vp,
 $$
 where $\cR_{|\nu|}$ is defined by \eqref{Eq4.5} and
 $$
\mathbf{P}^{\pm}(g)h= \frac{h \pm \boldsymbol\sigma(g)\,W(g_+^{-1})
\cR_{|\nu|}\, \vp}{2}.
 $$
 \end{thm}
 \textbf{Proof.}
Consider the expression $JQW^0(g)Ph$. One has
 \begin{align*}
JQW^0(g)Ph &=\boldsymbol\sigma(g) JQ  W^0(\widetilde{g}_+^{-1})
W^0(e^{i\nu t}) W^0 \left ( \left (\frac{t-i}{t+i}\right )^n\right )
W^0(g_+) Ph \\
 &=\boldsymbol\sigma(g) JQ  W^0(\widetilde{g}_+^{-1})
W^0(e^{i\nu t}) P\vp =\boldsymbol\sigma(g) P W^0(g_+^{-1})
W^0(e^{i|\nu| t}) J P\vp\\
 &=\boldsymbol\sigma(g) P W(g_+^{-1}) \cR_{|\nu|} \vp.
\end{align*}
Application of the relation
 $$
\mathbf{P}^{\pm}(g)h= \frac{h\pm JQW(g)P h}{2},
 $$
competes the proof.
  \rbx

 \section{Kernels and cokernels of Wiener-Hopf plus Hankel operators.
          Specification.\label{s3}}

In this section we study the kernels and cokernels of Wiener-Hopf
plus Hankel operators in the case where the generating functions
$a,b\in G$ satisfy the matching condition \eqref{cst9} and $a$ is
invertible in $G$. Then according to Theorem \ref{t1}, the operators
$W(c)$ and $W(d)$ are one-sided invertible in $L^p(\sR^+)$, $1\leq p
\leq \infty$. Using results of Section \ref{s2}, we derive an
explicit description for the kernels and cokernels of the operators
mentioned. As before, we again have to consider several cases.

 \subsection{The Case I: $\nu(c)=\nu(d)=0$.\label{ss2.1}}
This case is also used as a model case in order to show how to
handle all other situations. If the indexes $\nu(c)$ and $\nu(d)$
are equal to zero, then the operators $W(c)$ and $W(d)$ are
Fredholm. Using the relations (2.4) and (2.7) of \cite{DS:2014b},
one obtains that the operators $W(a)\pm H(b)$ are Fredholm.  Set
$\kappa_1:=\ind W(c)$, $\kappa_2:=\ind W(d)$ and let $\sZ_-$ and
$\sZ_+$ refer to the set of all negative and non-negative integers,
correspondingly.
  \begin{thm}\label{thm5}
 Assume that $\nu(c)=\nu(d)=0$.
 \begin{enumerate}
\item If $(\kappa_1,\kappa_2)\in \sZ_+\times \sN$, then for all  $p\in [1, \infty]$ the operators
$W(a)\pm H(b):L^p \to L^p$ are invertible from the right and
\begin{equation*}
\begin{aligned}
   \ker(W(a)+H(b))& =\im \mathbf{P}^-(c) \dotplus \vp_{+}(\im \mathbf{P}^+(d)), \\
\ker(W(a)-H(b))& =\im \mathbf{P}^+(c) \dotplus\vp_{-}(\im
\mathbf{P}^-(d)),
\end{aligned}
\end{equation*}
where the spaces $\im \mathbf{P}^{\pm}(c)$, $\im
\mathbf{P}^{\pm}(d)$ are described in Theorem \ref{thm1} and the
mappings $\vp_{\pm}$ are defined by \eqref{eqphi}.

\item If $(\kappa_1,\kappa_2)\in \sZ_-\times (\sZ\setminus \sN)$, then for all  $p\in [1, \infty]$ the operators
$W(a)\pm H(b):L^p \to L^p$ are invertible from the left and for all
$p\in [1, \infty)$ one has
\begin{equation*}
\begin{aligned}
   \coker(W(a)+H(b))& =\im \mathbf{P}^-(\overline{d}) \dotplus \vp_{+}(\im \mathbf{P}^+(\overline{c})), \\
\coker(W(a)-H(b))& =\im \mathbf{P}^+(\overline{d})
\dotplus\vp_{-}(\im \mathbf{P}^-(\overline{c})),
\end{aligned}
\end{equation*}
with $\im \mathbf{P}^{\pm}(\overline{d})=\{0\}$ if $\kappa_2=0$.

 \item If $(\kappa_1,\kappa_2)\in \sZ_+ \times (\sZ\setminus \sN)$, then for all  $p\in [1, \infty]$
one has
 \begin{equation*}
\begin{aligned}
   \ker(W(a)+H(b))& =\im \mathbf{P}^-(c), \\
\ker(W(a)-H(b))& =\im \mathbf{P}^+(c),
\end{aligned}
\end{equation*}
and for all $p\in [1, \infty)$,
\begin{equation*}
\begin{aligned}
   \coker(W(a)+H(b))& =\im \mathbf{P}^-(\overline{d}), \\
\coker(W(a)-H(b))& =\im \mathbf{P}^+(\overline{d}) ).
\end{aligned}
\end{equation*}
 \end{enumerate}
  \end{thm}
  \textbf{Proof.}  Let us note that all results concerning the kernels
of the corresponding operators follow immediately from Proposition
\ref{p2} and from Theorem \ref{thm1}. As far as the cokernel
structure is concerned, one has to take into account the already
mentioned relation \eqref{cst10} and the fact that
$(\overline{d},\overline{c})$ is the subordinated pair for the duo
$(\overline{a},\overline{b})$.
  \rbx

It remains to consider the case $(\kappa_1,\kappa_2)\in \sZ_- \times
\sN$. This situation is more involved. In order to formulate the
next result, we need a special representation for the index of the
operator $W(c)$. Thus chose $k\in \sN$ such that
 \begin{equation*} 1\geq
2k+\kappa_1\geq 0.
\end{equation*}
Such a number $k$ is uniquely defined and
\begin{equation*}
2k+\kappa_1 =\left\{%
\begin{array}{ll}
   0, & \hbox{if\;} \kappa_1 \; \hbox{is even,}  \\
    1, &\hbox{if\;} \kappa_1 \; \hbox{is odd.} \\
\end{array}%
\right.
\end{equation*}
Now the operators $W(a)\pm H(b)$ can be represented in the form
 \begin{equation}\label{eq6.1}
W(a)\pm H(b)= \left( W \left (a\left (\frac{t-i}{t+i} \right )^{-k}
\right )\pm H \left (b\left (\frac{t-i}{t+i} \right )^k \right
)\right )W \left (\left (\frac{t-i}{t+i} \right )^k \right ).
\end{equation}
Observe that $\left (a\left (\frac{t-i}{t+i} \right )^{-k}, b\left
(\frac{t-i}{t+i} \right )^k \right )$ is a matching pair with the
subordinated pair $\left (c\left (\frac{t-i}{t+i} \right )^{-2k},
d\right )$. Therefore, the operators $W\left (a\left
(\frac{t-i}{t+i} \right )^{-k}\right )\pm H \left (b\left
(\frac{t-i}{t+i} \right )^k \right )$ are subject to assertion (i)
of Theorem \ref{thm5}. Thus they are right-invertible, and if
$\kappa_1$ is even, then
 \begin{equation}\label{eq6.2}
\begin{aligned}
\ker \left (W\left (a\left (\frac{t-i}{t+i} \right )^{-k}\right)+
  H \left (b\left (\frac{t-i}{t+i} \right )^k\right)\right )=\vp_+(\im \mathbf{P}^{+}(d)),\\
 \ker \left (W\left (a\left (\frac{t-i}{t+i} \right )^{-k}\right)-
  H \left (b\left (\frac{t-i}{t+i} \right )^k\right) \right)=\vp_-(\im \mathbf{P}^{-}(d)),
\end{aligned}
\end{equation}
and if $\kappa_1$ is odd, then
\begin{equation}\label{eq6.3}
\begin{aligned}
\ker \left (W\left (a\left (\frac{t-i}{t+i} \right )^{-k}\right )
\right . & + \left .H \left (b\left (\frac{t-i}{t+i} \right
)^k \right) \right) \\
 &= \frac{1-\boldsymbol\sigma(c)}{2}W(c_+^{-1})\,\{\sC \psi_0\}
\dotplus  \vp_+(\im \mathbf{P}^{+}(d)),\\
\ker \left (W\left (a\left (\frac{t-i}{t+i} \right
)^{-k}\right)\right .&- \left. H \left (b\left (\frac{t-i}{t+i}
\right
)^k\right) \right) \\
&=\frac{1+\boldsymbol\sigma(c)}{2}W(c_+^{-1})\,\{\sC \psi_0\}
 \dotplus \vp_-(\im \mathbf{P}^{-}(d)),
\end{aligned}
\end{equation}
where the function $\psi_0$ is defined by \eqref{Eq4} and the
mappings $\vp_{\pm}$ depend on the functions $a\left
(\frac{t-i}{t+i} \right )^{-k}$ and $b\left (\frac{t-i}{t+i} \right
)^k$.

 \begin{thm}\label{t4}
Let $(\kappa_1,\kappa_2)\in \sZ_-\times \sN$ and $p\in [1,\infty)$.
Then
\begin{enumerate}
    \item If $\kappa_1$ is odd, then
   \begin{align*}
&\ker(W(a)\pm  H(b)) = W \left (\left (\frac{t-i}{t+i} \right
)^{-k}\right )\\
&\;\times\left (\left\{
 \frac{1\mp\boldsymbol\sigma(c)}{2}W(c_+^{-1})\,\{\sC \psi_0\}
\dotplus \vp_{\pm}(\im \mathbf{P}^{\pm}(d))\right\} \cap
\im W\! \left (\!\left (\frac{t-i}{t+i} \right )^k \right )\right )\\
 &\;=\left\{\!\!\psi\in \!\left \{ W \!\!\left (\left (\frac{t-i}{t+i} \right
)^{-k}\right )u\right \}\!: \!  u\in  \left
\{\frac{1\mp\boldsymbol\sigma(c)}{2}W(c_+^{-1})\,\{\sC \psi_0\}
 \dotplus \vp_{\pm}(\im \mathbf{P}^{\pm}(d))\right \} \right .  \\
 &\qquad\qquad \qquad\qquad\qquad\text{and}\left. \int_0^\infty u(t) e^{-t}t^j\,dt=0 \text{ for all } j=0,1,\cdots,
k-1,   \right \},
\end{align*}
where the mappings $\vp_{\pm}$ depend on the functions $a\left
(\frac{t-i}{t+i} \right )^{-k}$ and $b\left (\frac{t-i}{t+i} \right
)^k$. The last means that the functions $a,b$ and $c$ in the
expression \eqref{eqphi} have to be, respectively, replaced by
$a\left (\frac{t-i}{t+i} \right )^{-k}, b\left (\frac{t-i}{t+i}
\right )^k$ and $c\left (\frac{t-i}{t+i} \right )^{-2k}$.
    \item If $\kappa_1$ is even, then
     \begin{multline*}
 \ker(W(a)\pm  H(b))\!\! = \!\! W \left (\left (\frac{t-i}{t+i} \right )^{-k}\right)\!\!\!\!\left ( \left\{\vp_{\pm}(\im \mathbf{P}^{\pm}(d))\right\}
 \cap \im W \left (\left (\frac{t-i}{t+i} \right )^k\right ) \right )\\
 = \left\{ \psi\in \{ W \left (\left (\frac{t-i}{t+i} \right )^{-k}\right)u\}:
u\in \left \{\{\sC \psi_0\}\dotplus \vp_{\pm}(\im \mathbf{P}^{\pm}(d)) \, \right \}\right . \text{and} \\
  \left.
   \int_0^\infty u(t) e^{-t}t^j\,dt=0 \text{ for all } j=0,1,\cdots, k-1,   \right \},
 \end{multline*}
and the mappings $\vp_{\pm}$ again depend on $a\left
(\frac{t-i}{t+i} \right )^{-k}$ and $b\left (\frac{t-i}{t+i} \right
)^k$.
\end{enumerate}
\end{thm}

 \textbf{Proof.} 
It follows immediately from the representations
\eqref{eq6.1}--\eqref{eq6.3}.
\rbx

Theorem \ref{t4} can also be used to derive representations for the
cokernels of the operators $W(a)\pm H(b)$ in the situation where
$(\kappa_1,\kappa_2)\in \sZ_-\times \sN$. Indeed, recalling that for
$p\in [1,\infty)$, the adjoint operator  $(W(a)\pm H(b))^*$ can be
represented in the form \eqref{cst10} and $(\overline{d},
\overline{c})$ is the subordinated pair for
$(\overline{a},\widetilde{\overline{b}})$, one can observe that the
operators $W(\overline{d})$ and $W(\overline{c})$ are also Fredholm
and
\begin{equation*}
\ind W(\overline{d})=-\kappa_2, \quad \ind
W(\overline{c})=-\kappa_1,
\end{equation*}
so $(-\kappa_2,-\kappa_1)\in \sZ_-\times \sN$. Therefore, Theorem
\ref{t4} applies and one can formulate the following result.

 \begin{thm}\label{t5}
Let $(\kappa_1,\kappa_2)\in \sZ_-\times \sN$, and let $m\in\sN$
satisfy the requirement
\begin{equation*}
1\geq 2m-\kappa_2\geq0.
\end{equation*}
Then
\begin{enumerate}
    \item If $\kappa_2$ is odd, then
     \begin{align*}
& \coker(W(a)\pm H(b))= W \left (\!\left (\frac{t-i}{t+i} \right )^{-m}\right)\\
&\quad \times  \left (\left\{
 \frac{1\mp\boldsymbol\sigma(\overline{d})}{2}W(\overline{d_-^{-1}})\,
 \{\sC \psi_0\}
 \dotplus \vp_{\pm}(\im \mathbf{P}^{\pm}(\overline{c}))\right\} \cap \im W
 \left (\left (\frac{t-i}{t+i} \right )^{m}\right) \right ).
\end{align*}

    \item If $\kappa_2$ is even, then
     \begin{align*}
 &\coker(W(a)\pm H(b))= \\
 &= W\left (\left (\frac{t-i}{t+i} \right )^{-m}\right)\left (\left\{\vp_{\pm}(\im \mathbf{P}^{\pm}(\overline{c})\right\} \cap \im
 W\left (\left (\frac{t-i}{t+i} \right )^{m}\right) \right ),
 \end{align*}
and the mappings $\vp_{\pm}$ depend on $\overline{a}\left
(\frac{t-i}{t+i} \right )^{-m}$ and $\widetilde{\overline{b}}\left
(\frac{t-i}{t+i} \right )^{m}$.
\end{enumerate}
\end{thm}

\subsection{The Case II: $\nu(c)\neq 0$ and $\nu(d)\neq 0$.\label{ss2.2}}
According to Theorem \ref{t1}, the operators $W(c)$ and $W(d)$ are
one-sided invertible. In this situation the pair $(W(c),W(d))$
belongs to one of the classes $(r,r)$, $(l,l)$, $(l,r)$ or $(r,l)$,
where letter $r$ or $l$ means that the corresponding operator is
right- or left-invertible. It is worth mentioning that if the pair
$(W(c),W(d))$ belongs to the class $(r,l)$, then the operator
$W(a)+H(b)$ is normally solvable but it is not semi-Fredholm.
Further, if $(W(c),W(d))\in (l,r)$ then, generally, it is not known
whether $W(a)+H(b)$ is normally solvable or not. If $(W(c),W(d))$
belongs to one of the classes $(r,r)$ or $(r,l)$, then the kernels
of the operators $W(a)+H(b)$ and $W(a)-H(b)$ can be described using
results of Section \ref{s2}. For the description of the cokernels of
the operators $W(a)+H(b)$ and $W(a)-H(b)$ in the cases $(l,l)$ and
$(r,l)$, one has to assume that $p\in[1,\infty)$ and use the
relation \eqref{cst10}. If $(W(c),W(d))$ belongs to the class
$(r,l)$, then one can proceed similarly to Subsection \ref{ss2.1}.
More precisely, we have to consider three situations, namely,
\begin{enumerate}
    \item  The index $\nu(c)<0$ and the index $n(c)>0$.
    \item The index $\nu(c)<0$ and the index $n(c)=0$.
    \item The index $\nu(c)<0$ and the index $n(c)<0$.
\end{enumerate}
Since in this situations, the operator $W(c)$ is right-invertible,
the kernels of the operators $W(a)+H(b)$ and $W(a)-H(b)$ can be
described by Proposition \ref{p2} and subsequent use of Theorems
\ref{thm2},  \ref{thm3} and \ref{thm4}.

As was already mentioned, if the pair $(W(c),W(d))$ belongs to the
class $(l,r)$, then it is not known whether the operators $W(a)\pm
H(b)$ are normally solvable. Nevertheless, the kernels and cokernels
of these operators still can be described. However, it is worth
noting that Proposition \ref{p2} cannot be directly used. Thus let
us sketch the idea how to proceed in this situation. We have to deal
with the following cases
  \begin{enumerate}
    \item  The index $\nu(c)>0$ and the index $n(c)>0$.
    \item The index $\nu(c)>0$ and the index $n(c)=0$.
    \item The index $\nu(c)>0$ and the index $n(c)<0$.
\end{enumerate}
In these situations the operators $W(a)\pm H(b)$ admit the
factorization
  \begin{align*}
W(a)\pm H(b) &= \left ( W \left ( ae^{-i\nu t/2} \left (
\frac{t-i}{t+i} \right)^{-k} \right )\pm  H\left ( be^{i\nu t/2}
\left ( \frac{t-i}{t+i} \right)^{k} \right )  \right ) \\ &
\qquad\qquad \times W \left (
e^{i\nu t/2} \left ( \frac{t-i}{t+i} \right)^{k} \right ),\\
W(a)\pm H(b) &= \left ( W \left ( ae^{-i\nu t/2} \right )\pm H\left
(
 be^{i\nu t/2} \right )  \right ) W \left ( e^{i\nu t/2} \right
 ), \\
 W(a)\pm H(b) &= \left ( W \left ( ae^{-i\nu t/2} \right )\pm  H\left (
 be^{i\nu t/2} \right )  \right )  W \left ( e^{i\nu t/2} \right
 ),
\end{align*}
where $\nu=\nu(c)$ and $k$ are defined as in Subsection \ref{ss2.1}.
Let us consider, for definiteness, the operator $W(a)+H(b)$ and note
that the respective subordinated pairs for the first operators in
the right-hand sides of the last representations are
 \begin{equation*}
  \left ( ce^{-i\nu t} \left (
\frac{t-i}{t+i} \right)^{-2k} , d \right ), \left ( ce^{-i\nu t}, d
\right), \text{ and } \left ( ce^{-i\nu t}, d \right)
 \end{equation*}
with the respective indices $\nu$ and $n$ defined as
 \begin{align*}
& \nu  \left ( ce^{-i\nu t} \left ( \frac{t-i}{t+i} \right)^{-2k}
\right ) =0 \text { and } n\left ( c
e^{-i\nu t} \left ( \frac{t-i}{t+i}\right )^{-2k}\right ) =-2k+n(c), \\
& \nu  \left ( ce^{-i\nu t} \right ) =0 \text { and } n\left ( c
e^{-i\nu t} \right ) =0,\\
 & \nu  \left ( ce^{-i\nu t} \right ) =0 \text { and } n\left ( c
e^{-i\nu t} \right ) =n(c).
\end{align*}
Now using the corresponding results of Section \ref{s2} and those
obtained in Subsection \ref{ss2.1}, one can get a complete
description for the kernels and cokernels of the operators
$W(a)+H(b)$ and $W(a)-H(b)$.

 \subsection{The Case III.\label{ss2.3}}
Assume that the only one of the indices $\nu(c)$ or $\nu(d)$ is
equal to zero. This case can be handled similarly to the Cases~I and
II without any new features. Therefore, we omit detailed
formulations here. However, it is worth mentioning that in this case
the operators $W(a)\pm H(b)$ are semi-Fredholm but not Fredholm.

\section*{Acknowledgements}

The authors thank an anonymous referee for very careful reading of
the manuscript and suggesting a number of improvements and
corrections.


 \end{document}